\documentclass[11pt,bibay3]{asl}

\UseRawInputEncoding
\makeatletter

\def\bottomnote{\@ifnextchar
    [{\@xfootnotenext}{\xdef\@thefnmark{}\@footnotetext}}
\def\@xfootnotenext[#1]{\begingroup \csname c@\@mpfn\endcsname #1\relax
   \xdef\@thefnmark{\thempfn}\endgroup \@footnotetext}

\newcount\currhour
\newcount\currminutes
\newcount\temp
\currhour=\time\currminutes=\time
\divide\currhour by 60 
\temp=\currhour
\multiply\temp by 60    
\advance\currminutes by -\temp
\def\minutes{\ifnum\currminutes<10 0\number\currminutes%
\else\number\currminutes\fi}
\def\now{\today, \number\currhour:\minutes}

\def\plaindraft#1{%
\gdef\draftinfo{#1}
\let\oldps@headings\ps@headings%
\def\ps@headings{\oldps@headings%
\def\@oddfoot{\vbox to 0pt{\vskip 10pt\noindent \small\tt
#1 }}%
\def\@evenfoot{\vbox to 0pt{\vskip 10pt\noindent \small\tt
#1 }}}%
\AtBeginDocument{\ttdraft}%
\pagestyle{headings}}

\def\myplaindraft#1{%
\gdef\draftinfo{#1}
\let\oldps@myheadings\ps@myheadings%
\def\ps@myheadings{\oldps@myheadings%
\def\@oddfoot{\vbox to 0pt{\vskip 10pt\noindent \small\tt
#1 }}%
\def\@evenfoot{\vbox to 0pt{\vskip 10pt\noindent \small\tt
#1 }}}
\pagestyle{myheadings}}

\def\draft#1{\plaindraft{#1 \now\hfill\thepage}}

\def\ttdraft{%
\let\oldps@firstpage\ps@firstpage%
\gdef\ps@firstpage{\oldps@firstpage
\def\@oddfoot{\vbox to 0pt{\vskip 10pt\noindent
\small\tt \draftinfo}}}}

\newcommand\yf{\hspace*{\fill}}

\newcommand\Graph{\textup{Graph}}

\def\draft#1{%
\gdef\draftinfo{#1}
\let\oldps@headings\ps@headings%
\def\ps@headings{\oldps@headings%
\def\@oddfoot{\vbox to 0pt{\vskip 10pt\noindent \small\tt
\draftinfo \hfill\thepage}}%
\def\@evenfoot{\vbox to 0pt{\vskip 10pt\noindent \small\tt
\draftinfo \hfill\thepage}}}
\let\oldps@firstpage\ps@firstpage%
\gdef\ps@firstpage{\oldps@firstpage
\def\@oddfoot{\vbox to 0pt{\vskip 10pt\noindent
\small\tt \draftinfo \hfill\thepage}}}
\pagestyle{headings}}    

\makeatother

\newcommand{\ccX}{{\mathcal{X}}}
\newcommand{\ccY}{{\mathcal{Y}}}

\newcommand\cf{cf.~}
\newcommand{\lh}{\mbox{\rm lh}}
\newcommand{\baire}{\mathcal{N}}
\newcommand{\eedf}{\iff\kern -.65em_{\mathrm{df}}~}
\newcommand\st{ : }
\newcommand\vep{\varepsilon}
\newcommand\tb[1]{\tboldsymbol{#1}}
\newcommand\nats{\mathbb{N}}
\newcommand{\bij}{\mbox{\,$\rightarrowtail\kern -.8em \rightarrow$\,}}
\newcommand{\pfto}{\rightharpoonup}
\def\seq#1{\langle {#1}\rangle}
\newcommand\Seq{\mathrm{Seq}}
\newcommand{\converges}{\kern -.2em\downarrow\kern .1em}
\newcommand{\diverges}{\kern-.1em\uparrow\kern .05em}
\newcommand{\Det}{{\textup{Det}}}
\newcommand\conj{{~\&~}}
\newcommand\ykeep[1]{\makebox[0pt][l]{#1}}
\newcommand\oo[1]{\overline{#1}}

\title[The strong Spector-Gandy Theorem]
{The strong Spector-Gandy Theorem for the higher analytical pointclasses${\,}^\ast$}

\author{Joan R. Moschovakis}
\address{Occidental College}
\email{joan.rand@gmail.com}
\urladdr{http://www.math.ucla.edu/\urltilde joan}

\author{Yiannis N. Moschovakis} \twoaddress{
Department of Mathematics\\
University of California,
Los Angeles, CA, USA}
{Department of Mathematics,
University of Athens, Greece}{3}
\email{ynm@math.ucla.edu}
\urladdr{http://www.math.ucla.edu/\urltilde ynm}

\plaindraft{sg, Feb. 5, 2022, \hfill\now p. \thepage}

\def\yextra{}
\newcommand\dst[1]{\textup{DST{#1}}}

\newcommand\ps{p}
\newcommand\bu{\mathbf{u}}
\newcommand\bv{\mathbf{v}}

\newtheorem*{klrestr}{\bfseries Restricted Quantification Theorem}
\newtheorem*{sgthm}{(SG)~\bfseries Spector-Gandy Theorem}
\newtheorem*{ssgthm}{(SSG)~\bfseries Strong Spector-Gandy Theorem}
\newtheorem*{extssg1}{\bfseries Strong Spector-Gandy Theorem for $\Pi^1_{2n+1}$}
\newtheorem*{extssg2}{\bfseries Main Theorem}
\newtheorem*{third}{\bfseries Third Periodicity Theorem}

\newtheorem{theorem}{\bfseries Theorem\yextra}

\theoremstyle{definition}

\begin{document}

\maketitle


\bottomnote{${}^\ast$The results in this note will be included in
Part III of \cite{edst}, the advanced part for which we assume
familiarity with the material (and notation) of Chapters 1~--~7
of~\cite{dst} (\dst).}

There is little doubt that the single, most important elementary
result of effective descriptive set theory is the

\begin{klrestr}[\cite{kleene1959b}, \dst{(4D.3)}].\\
For any two spaces $\ccX,
\ccY$ of type $\leq1$, if\/ $Q\subseteq\ccX\times\ccY$ is $\Pi^1_1$
and\footnote{A space of type $\leq1$ is a finite product of copies of $\nats$
and $\baire$.}
\[
P(x) \iff (\exists y\in\Delta^1_1(x))Q(x,y),
\]
then $P$ is also $\Pi^1_1$.
\end{klrestr}

It yields the $\Delta$-\textit{Uniformization Criterion}
\dst{(4D4)}, which is the key tool in the effective proofs of
several results about Borel uniformizations  of
Borel sets, \cf Section 4F of \dst. It also motivated Spector to
prove the following basic fact which provides a converse
to it:

\begin{sgthm}[\cite{spector1960}, also \cite{gandy1960}]
For every space $\ccX$ of type $\leq1$ and every $\Pi^1_1$ set $P\subseteq\ccX$,
there is some $\Pi^0_1$ set $Q\subseteq\ccX\times\baire$ such that
\begin{equation}
\label{ast2}
P(x) \iff (\exists\alpha\in\Delta^1_1(x))Q(x,\alpha);
\end{equation}
and so
\begin{equation}
\label{ast3}
P\in\Pi^1_1(\ccX) \iff \text{\eqref{ast2} holds with some }
Q\in\Pi^0_1(\ccX\times\baire).
\end{equation}
\end{sgthm}

Spector actually proved the following stronger form of this
theorem:

\begin{ssgthm}[\cite{spector1960}]{}
For every space $\ccX$ of type $\leq1$  and every $\Pi^1_1$ set $P\subseteq\ccX$,
there is some $\Pi^0_1$ set $Q\subseteq\ccX\times\baire$ such that
\begin{equation}
\label{ssgeq}
P(x) \iff (\exists!\alpha)Q(x,\alpha)\iff (\exists\alpha\in\Delta^1_1(x))Q(x,\alpha).
\end{equation}
\end{ssgthm}

In \dst{(6E.7)}, ynm gave a proof of (SG) for the odd levels
$\Pi^1_{2n+1}$ of the analytical hierarchy under the hypothesis of
\textit{Projective Determinacy} and our aim here is to prove a
similar extension of (SSG) as follows:

\begin{extssg1}
\label{mainthm}
Assume $\Det(\tb\Sigma^1_{2n})$. For every space $\ccX$ of type
$\leq1$ and every $\Pi^1_{2n+1}$ set $P\subseteq\ccX$, there is
a $\Pi^1_{2n}$ set $Q\subseteq\ccX\times\baire$ such that
\begin{equation}
\label{mainthmeq}
P(x)\iff(\exists!\alpha)Q(x,\alpha)
\iff(\exists\alpha\in\Delta^1_{2n+1}(x))
Q(x,\alpha).
\end{equation}
\end{extssg1}

There were many reasons why Gandy's weaker  result  (without the
$(\exists!\alpha)$) was generally adopted in the 70's as the
standard version of this basic fact: it is quite easy to prove; it
suffices to prove~\eqref{ast3}; and (most significantly), it can
be established in a similar way for the analog of $\Pi^1_1$ in many theories of
effective definability, including \textit{normal} (Kleene)
\textit{recursion} in type $3$ (Theorem~9 in \cite{ynm1967}) and
\textit{inductive definability} in every acceptable structure
$\mathfrak{A}$, Theorem 7D.2 in \cite{ynmeias}.
What was missed in the 70's is that (SSG) has important
applications in effective descriptive set theory, \cf Chapter~5 of
\cite{gregoriades2016}; it remains to be seen if the same is true
of its extension to the higher analytical pointclasses under
projective determinacy.

\section{Background}
\label{background}
The proof outlined in the next section uses the methods and results
of Sections~6D and~6E of DST about the game quantifier
$\Game{\alpha}$ defined in the beginning of 6D,
\begin{multline*}
(\Game \alpha)Q(x,\alpha) \eedf \text{ Player $I$ wins the game }\{\alpha\st Q(x,\alpha)\}\\
\iff (\exists a_0)(\forall a_1)(\exists a_2) \cdots Q(x,(a_0, a_1, a_2, \ldots)).
\end{multline*}
In particular, we will use systematically the (trivial) characterizations
\[
\Pi^1_1 = \Game\Sigma^0_1,~~\Sigma^1_2=\Game\Pi^1_1, ~~\Pi^1_3=\Game \Sigma^1_2, \ldots
\]
rather than the classical definitions of these pointclasses in
terms of the quantifiers $\exists\alpha, \forall\alpha$. We will
work only with spaces $\ccX$ of type $\leq1$, we fix \textit{good
parametrizations in $\baire$} (\cf \dst{(3H.1)}) in some
standard way on $\Sigma^1_0=\Sigma^0_1$ and we extend them to all
the analytical pointclasses by the recursion
\begin{gather*}
G^{\Pi,1}_{2n+1}(\vep,x) \eedf (\Game\alpha)G^{\Sigma,1}_{2n}(\vep,x,\alpha),
\\
G^{\Sigma,1}_{2n+2}(\vep,x) \eedf (\Game\alpha)G^{\Pi,1}_{2n+1}(\vep,x,\alpha),
\\
G^{\Sigma,1}_{2n+1}(\vep,x) \eedf \lnot G^{\Pi,1}_{2n+1}(\vep,x),
\quad
G^{\Pi,1}_{2n}(\vep,x) \eedf \lnot G^{\Sigma,1}_{2n}(\vep,x);
\end{gather*}
$\vep$ is a $\tb\Delta^1_k$ code of some $Q\subseteq\ccX$ if
\[
Q(x) \iff G^{\Sigma,1}_k((\vep)_0,x) \iff G^{\Pi,1}_k((\vep)_1,x).
\]

The first, basic theorem we need is (perhaps) the main fact about the game quantifier:

\begin{third}[\dst{(6E.1)}, ynm]
\label{third}
Assume $\Det(\tb\Sigma^1_{2n})$. If a pointset $A\subseteq\baire$ is in $\tb\Sigma^1_{2n}$
with code $\zeta$ and player I wins $A$, then I has a $\Delta^1_{2n+1}(\zeta)$ winning
strategy.
\end{third}

The proof of \dst{(6E.1)} gives a great deal of information about
$\tb\Pi^1_{2n+1}$ beyond its statement, and  we need to make
explicit and use some of this here.

A \textit{partial strategy} (for player $I$ in a game on $\nats$)
is a partial function $\ps : \nats\pfto\nats$ such that
$\lnot(\Seq(u)\conj\lh(h)\text{ is even})\implies \ps(u)=0$; if $\ps$ is total,
it then defines a strategy for player $I$ by which he plays
\[
a_0=\ps(\seq{~}),~a_2=\ps(\seq{a_0, a_1}), a_4 = \ps(\seq{a_0, a_1,a_2, a_3}) \ldots
\]
when Player $II$ plays  $a_1, a_3, \ldots$.
As usual, for any $p:\nats\pfto\nats$,
\[
\Graph(p) = \{(u,t) \st p(u)=t\}
\]
and we typically identify $p$ with its graph.

\begin{theorem}
\label{fact1}

If\/ $\Det(\tb\Sigma^1_{2n})$ and $\ccX$ is of type $\leq 1$, then
there is a recursive function
\[
\bu_1(\vep,x)=\seq{\bu_1^\Sigma(\vep,x),\bu_1^\Pi(\vep,x)}
\]
and for each $\vep$ and each $x$, a partial function $\ps_{\vep,x}$ such that
\begin{equation}
\label{fact1eq1}
\ps_{\vep,x}(u)=t \iff G^{\Pi,1}_{2n+1}(\bu_1^\Pi(\vep,x),u,t),
\end{equation}
and if \qquad${
P_\vep(x) \eedf G^{\Pi,1}_{2n+1}(\vep,x)
\iff (\Game\alpha)G^{\Sigma,1}_{2n}(\vep,x,\alpha)}$,\qquad then
\begin{align}
\label{fact1eq2}
\quad P_\vep(x) &\iff\ps_{\vep,x}\text{ is total }
\iff(\exists\beta)(\forall u)G^{\Pi,1}_{2n+1}(\bu_1^\Pi(\vep,x),u,\beta(u))\\
\label{fact1eq3}
&\iff (\exists!\beta)(\forall u)G^{\Pi,1}_{2n+1}(\bu_1^\Pi(\vep,x),u,\beta(u))\\
\label{fact1eq4}
&\iff
\Graph(\ps_{\vep,x}) \text{ is $\tb\Delta^1_{2n+1}$ with code }\bu_1(\vep,x)\\
\label{fact1eq6}
&\iff (\exists\beta\in\Delta^1_{2n+1}(\vep,x))(\forall u)G^{\Pi,1}_{2n+1}(\bu_1^\Pi(\vep,x),u,\beta(u)).
\end{align}
\end{theorem}

The construction  of $\ps_{\vep,x}$ is quite complex but
effective, and so the definition of $\bu^\Pi_1(\vep,x)$ so
that~\eqref{fact1eq1} holds is routine. The same construction
proves~\eqref{fact1eq2}, which is trivially equivalent
with~\eqref{fact1eq3}. To define $\bu_1^\Sigma(\vep,x)$ and prove
the crucial~\eqref{fact1eq4}, we use the closure of $\Sigma^1_{2n+1}$
under number quantification and the characteristic property of a good
parametrization, which supplies $\bu^\Sigma_1(\vep,x)$ so that
for all $\vep,x,u,t$,
\begin{equation}
\label{goodparam}
G^{\Sigma,1}_{2n+1}(\bu^\Sigma_1(\vep,x),u,t) \iff (\forall s)
\Big(G^{\Pi,1}_{2n+1}(\bu_1^\Pi(\vep,x),u,s)\implies t = s\Big),
\end{equation}
and this implies easily that
\begin{multline}
\label{fact1eq5}
P_\vep(x) \implies \ps_{\vep,x}\text{ is total }\\
\implies (\forall u)(\forall t)\Big(
G^{\Sigma,1}_{2n+1}(\bu_1^\Sigma(\vep,x),u,t)
\iff G^{\Pi,1}_{2n+1}(\bu_1^\Pi(\vep,x),u,t)\Big).
\end{multline}
To prove the converse of~\eqref{fact1eq5} and complete the proof of~\eqref{fact1eq4},
assume the right-hand-side holds, but
(towards a contradiction), $\ps_{\vep,x}$ is not total,
so for some $\oo u$, $\ps_{\vep,x}(\oo u)\diverges$; but then
\[
(\forall s)\lnot G^{\Pi,1}_{2n+1}(\bu_1^\Pi(\vep,x),\oo u,s)
\text{ but }G^{\Sigma,1}_{2n+1}(\bu_1^\Sigma(\vep,x),\oo u,0),
\]
which is absurd.

Finally,~\eqref{fact1eq4} with~\eqref{fact1eq2} trivially imply~\eqref{fact1eq6}.\smallskip

The second theorem we need is a basic fact about the structure of
the analytical hierarchy under determinacy hypotheses:

\begin{theorem}[\dst{(6E.14)}, ynm]
\label{fact2}

If\/ $\Det(\tb\Sigma^1_{2n})$ and $\ccX$ is
of type $\leq 1$, then there is a recursive function
$\bu_2(\zeta)$, such that for all \ykeep{$R\subseteq\ccX$ and all $\zeta$:}
\begin{multline*}
\text{if $R$ is $\tb\Delta^1_{2n+1}$ with code $\zeta$, then}\\
\qquad R(x) \iff (\exists\gamma)G^{\Pi,1}_{2n}(\bu_2(\zeta), x,\gamma) \iff
(\exists!\gamma)G^{\Pi,1}_{2n}(\bu_2(\zeta), x,\gamma)\hspace*{\fill}\\
\hspace*{\fill}\iff(\exists\gamma\in\Delta^1_{2n+1}(\zeta,x))G^{\Pi,1}_{2n}(\bu_2(\zeta), x,\gamma).
\end{multline*}
\end{theorem}

Theorem~\ref{fact2} is an effective version of the
extension to $\tb\Delta^1_{2n+1}$ of Lusin's favorite
characterization of the Borel sets, \dst{(2E.8)}:
\begin{multline*}
R\subseteq\ccX \textit{ is Borel}\\
\iff R \textit{ is the continuous, injective
image of a closed }F\subseteq\baire.
\end{multline*}
\noindent
Its proof uses most everything preceding it in Section 6 of \dst,
including the version of the 3rd Periodicity Theorem \dst{(6E.1)}
in Theorem~\ref{fact1} and \textit{Solovay's trick}, as  in the
proof of \dst{(6E.12)}. However, the computation of codes we need
for this uniform version is quite routine.

\section{The main result}

We assume the notation of \S 1 and prove the following, uniform version of the Strong
Spector-Gandy Theorem for $\Pi^1_{2n+1}$ on page~\pageref{mainthm}:

\begin{extssg2}
\label{fact3}

If\/ $\Det(\tb\Sigma^1_{2n})$ and $\ccX$ is of type
$\leq1$, then there is a recursive function $\bv(\vep)$
such that for all $\vep$ and all $x$,
\begin{multline}
\label{finalversion}
G^{\Pi,1}_{2n+1}(\vep,x) \iff (\exists!\alpha)G^{\Pi,1}_{2n}(\bv(\vep),x,\alpha)\\
 \iff
(\exists\alpha\in\Delta^1_{2n+1}(\vep,x))G^{\Pi,1}_{2n}(\bv(\vep),x,\alpha).
\end{multline}
\end{extssg2}

Crucial to the proof will be the following equivalence which holds
for every $R\subseteq\baire\times\baire$ (jrm):
\begin{equation}
\label{joan}
(\exists!\beta)(\exists\gamma)R(\beta,\gamma)\conj(\exists\beta)(\exists!\gamma)R(\beta,\gamma)
\iff (\exists!(\beta,\gamma))R(\beta,\gamma),
\end{equation}
where the \textit{binary quantifier} on the right is defined in the natural way,
\begin{multline*}
(\exists!(\beta,\gamma))R(\beta,\gamma) \eedf (\exists\beta)(\exists\gamma)R(\beta,\gamma)\\
\conj (\forall\beta)(\forall\gamma)(\forall \beta')(\forall \gamma')\Big
([R(\beta,\gamma)\conj R(\beta',\gamma')]\implies [\beta=\beta'\conj\gamma=\gamma']\Big).
\end{multline*}

\textit{Proof of}~\eqref{joan}.
$(\implies)$ Let $\oo\beta$ be the unique $\beta$ such that $(\exists\gamma)R(\oo\beta,\gamma)$;
if $\beta'$ is such that $(\exists!\gamma)R(\beta',\gamma)$, then $(\exists\gamma)R(\beta',\gamma)$,
so $\beta'=\oo\beta$. So the left-hand-side gives unique $\oo\beta$ and $\oo\gamma$ for which
$R(\oo\beta,\oo\gamma)$ and then the right-hand-side holds with the pair $(\oo\beta,\oo\gamma)$.

$(\Longleftarrow)$ Let $(\oo\beta,\oo\gamma)$ be the unique pair
$(\beta,\gamma)$ such that $R(\oo\beta,\oo\gamma)$. This implies that
$\oo\beta$ is the unique $\beta$ such that for some $\gamma$,
$R(\oo\beta,\gamma)$, which implies that there is such a $\beta$,
which gives the first conjunct in the left-hand-side. In the same
way, there is a $\beta$ (namely $\oo\beta$) such that
$(\exists!\gamma)R(\oo\beta,\gamma)$, which proves the second
conjunct on the left. \yf\qedsymbol~\eqref{joan}\smallskip

To interpret the binary quantifier $(\exists!(\beta,\gamma))$ in terms of
the unary $(\exists!\alpha)$, we need a recursive isomorphism of $\baire$
with $\baire\times\baire$, and the simplest of these is the ``interweaving'' of two sequences
\[
\alpha\mapsto ((\alpha(0), \alpha(2), \ldots, ),(\alpha(1), \alpha(3), \ldots))
= (\lambda(t)\alpha(2t), \lambda(t)\alpha(2t+1));
\]
using this,~\eqref{joan} takes the form
\begin{multline}
\label{joan2}
(\exists!\beta)(\exists\gamma)R(\beta,\gamma)\conj(\exists\beta)(\exists!\gamma)R(\beta,\gamma)
\\\iff (\exists!\alpha)R(\lambda(t)\alpha(2t),\lambda(t)\alpha(2t+1)),
\end{multline}
which is the form of this equivalence we will use in the proof that follows.

\begin{proofplain} of the Main Theorem.
Fix $\vep$, let $P_\vep(x) \eedf G^{\Pi,1}_{2n+1}(\vep,x)$ and set
$\bu(\vep,x) = \bu_2(\bu_1(\vep,x))$.
By~\eqref{fact1eq4} and~\eqref{fact1eq1} in Theorem~\ref{fact1},
\begin{multline*}
P_\vep(x) \iff \{(u,t) \st G^{\Pi,1}_{2n+1}(\bu_1^\Pi(\vep,x),u,t)\}\\
\text{is $\tb\Delta^1_{2n+1}$ with code }\bu_1(\vep,x);
\end{multline*}
and then by Theorem~\ref{fact2} with $\zeta:=\bu_1(\vep,x)$ and for all $u$ and $t$,
\begin{multline}
\label{alt1}
G^{\Pi,1}_{2n+1}(\bu_1^\Pi(\vep,x),u,t) \iff (\exists\gamma)G^{\Pi,1}_{2n}(\bu(\vep,x),u,t,\gamma)\\
\yf\iff (\exists!\gamma)G^{\Pi,1}_{2n}(\bu(\vep,x),u,t,\gamma);
\end{multline}
and so for all $u$ and $\beta$,
\begin{multline}
\label{alt2}
G^{\Pi,1}_{2n+1}(\bu_1^\Pi(\vep,x),u,\beta(u)) \iff (\exists\gamma)G^{\Pi,1}_{2n}(\bu(\vep,x),u,\beta(u),\gamma)\\
\yf\iff (\exists!\gamma)G^{\Pi,1}_{2n}(\bu(\vep,x),u,\beta(u),\gamma).
\end{multline}
Using the first equivalence in~\eqref{alt2}  with~\eqref{fact1eq3} in Theorem~\ref{fact1}, we get
\begin{multline}
\label{alt3}
P_\vep(x) \iff (\exists!\beta)(\forall u)(\exists\gamma) G^{\Pi,1}_{2n}(\bu(\vep,x),u,\beta(u),\gamma)\\
\iff (\exists!\beta)(\exists\gamma)(\forall u)G^{\Pi,1}_{2n}(\bu(\vep,x),u,\beta(u),(\gamma)_u);
\end{multline}
and using~\eqref{fact1eq2} of Theorem~\ref{fact1} with the second one, we get
\begin{multline}
\label{alt4}
P_\vep(x) \iff (\exists\beta)(\forall u)(\exists!\gamma) G^{\Pi,1}_{2n}(\bu(\vep,x),u,\beta(u),\gamma)\\
\iff (\exists\beta)(\exists!\gamma)(\forall u)G^{\Pi,1}_{2n}(\bu(\vep,x),u,\beta(u),(\gamma)_u).
\end{multline}
If we now set
\begin{equation}
\label{R}
R_\vep(x,\beta,\gamma) \eedf (\forall u)G^{\Pi,1}_{2n}(\bu(\vep,x),u,\beta(u),(\gamma)_u)
\end{equation}
and put~\eqref{alt3} and~\eqref{alt4} together, we get
\begin{equation}
\label{alt5}
P_\vep(x) \iff (\exists!\beta)(\exists\gamma)R_\vep(x,\beta,\gamma)
\conj (\exists\beta)(\exists!\gamma)R_\vep(x,\beta,\gamma);
\end{equation}
and then~\eqref{joan2} gives
\begin{multline*}
P_\vep(x) \iff (\exists!\alpha)R_\vep(x,\lambda(t)\alpha(2t),\lambda(t)\alpha(2t+1))
\\\iff (\exists!\alpha)\Big[(\forall u)G^{\Pi,1}_{2n}(\bu(\vep,x),u,\lambda(t)\alpha(2t)(u),
(\lambda(t)\alpha(2t+1))_u)\Big].
\end{multline*}
Finally, the relation within the braces $\Big[~~\Big]$
is $\Pi^1_{2n}$, and so the good parametrization
property  gives us a recursive $\bv(\vep)$ such that
\begin{equation}
\label{bv}
P_\vep(x) \iff (\exists!\alpha)G^{\Pi,1}_{2n}(\bv(\vep),x,\alpha),
\end{equation}
which is the first claim in the theorem.\smallskip

The second claim is proved in the same way, with a couple of the
steps justified by different arguments.\smallskip

By~\eqref{fact1eq4} and~\eqref{fact1eq1} in Theorem~\ref{fact1},
\begin{multline*}
P_\vep(x) \iff \{(u,t) \st G^{\Pi,1}_{2n+1}(\bu_1^\Pi(\vep,x),u,t)\}\\
\text{is $\tb\Delta^1_{2n+1}$ with code }\bu_1(\vep,x);
\end{multline*}
and then by Theorem~\ref{fact2}, setting $\zeta:=\bu_1(\vep,x)$, for all $u$ and $t$,
\begin{multline*}
G^{\Pi,1}_{2n+1}(\bu_1^\Pi(\vep,x),u,t) \implies
(\exists\gamma\in\Delta^1_{2n+1}(\bu_1(\vep,x)  ))G^{\Pi,1}_{2n}(\bu(\vep,x),u,t,\gamma)\\
\implies \underline{(\exists\gamma\in\Delta^1_{2n+1}(\vep,x))G^{\Pi,1}_{2n}(\bu(\vep,x),u,t,\gamma)}\\
\qquad\implies (\exists\gamma)G^{\Pi,1}_{2n}(\bu(\vep,x),u,t,\gamma)\\
\implies G^{\Pi,1}_{2n+1}(\bu_1^\Pi(\vep,x),u,t),
\end{multline*}
where the step introducing the underlined proposition is justified because, as a general
fact,
\[
\Delta^1_{2n+1}(\bu_1(\vep,x))\subseteq \Delta^1_{2n+1}(\vep,x).
\]
It follows that for all $u$ and $t$,
\[
G^{\Pi,1}_{2n+1}(\bu_1^\Pi(\vep,x),u,t) \iff
(\exists\gamma\in\Delta^1_{2n+1}(\vep,x))G^{\Pi,1}_{2n}(\bu(\vep,x),u,t,\gamma),
\]
and hence for all $u$ and $\beta$,
\begin{multline*}
G^{\Pi,1}_{2n+1}(\bu_1^\Pi(\vep,x),u,\beta(u))\\ \iff
(\exists\gamma\in\Delta^1_{2n+1}(\vep,x))G^{\Pi,1}_{2n}(\bu(\vep,x),u,\beta(u),\gamma).
\end{multline*}

Using this with~\eqref{fact1eq6} in Theorem~\ref{fact1}, we get
\begin{multline*}
P_\vep(x) \iff
(\exists\beta\in\Delta^1_{2n+1}(\vep,x))(\forall u)
(\exists\gamma\in\Delta^1_{2n+1}(\vep,x))\\
G^{\Pi,1}_{2n}(\bu(\vep,x),u,\beta(u),\gamma);
\end{multline*}
and from this, it is a standard exercise in restricted quantification to get
\begin{multline*}
P_\vep(x) \iff
(\exists\beta\in\Delta^1_{2n+1}(\vep,x))
(\exists\gamma\in\Delta^1_{2n+1}(\vep,x))\\
\yf(\forall u)G^{\Pi,1}_{2n}(\bu(\vep,x),u,\beta(u),(\gamma)_u)\\
\iff (\exists\alpha\in\Delta^1_{2n+1}(\vep,x))
(\forall u)G^{\Pi,1}_{2n}(\bu(\vep,x),u,\lambda(t)\alpha(2t)(u),(\lambda(t)\alpha(2t+1))_u)\\
\iff (\exists\alpha\in\Delta^1_{2n+1}(\vep,x))G^{\Pi,1}_{2n}(\bv(\vep),x,\alpha)
\end{multline*}
where $\bv(\vep)$ is introduced in~\eqref{bv}, and this is the second claim in the theorem.
\end{proofplain}

\section{Comments}

The cumbersome computation of ``uniformities'' like $\bu_1(\vep,x)$,
$\bu_2(\zeta), \ldots, \bv(\vep)$ is rarely needed to discover,
prove  and communicate results in effective descriptive set
theory.

Consider, for example, a proposition of the form
\begin{equation}
\tag{$\ast$}
\label{star1}
(\forall x)(\exists y\in\Delta^1_1(x)) P(x,y)\qquad (P\subseteq\ccX\times\ccY).
\end{equation}

When one is primarily interested in applications, one replaces
it by the weaker
\begin{equation}
\tag{$\ast\ast$}
\label{star2}
(\exists \text{ a Borel function } f:\ccX\to\ccY)(\forall x)P(x,f(x))
\end{equation}
which (typically) suffices for the application we are considering.
A proof of~\eqref{star1} often yields a \textit{classical
proof} of~\eqref{star2} which can be  understood without direct
reference to the effective notions. For a good example of this,
see Corollary (39.20) in \cite{kechris1995}, which is a
``classical version'' of the Third Periodicity Theorem
\dst{(6E.1)} as we quoted it on page~\pageref{third}.

In the effective theory, where we want to keep track of the
parameters involved, the key, intuitive notion is that of a
\textit{constructive} (or \textit{effective}, or \textit{uniform})
\textit{proof}. For example, to prove the first claim in the Main
Theorem, one tries to prove \textit{constructively} that
\begin{multline*}
P(x) \text{ is }\tb\Pi^1_{2n+1}
\iff \textit{ there is a set $R(x,\alpha)$ in $\tb\Pi^1_{2n}$ such that }\\
P(x) \iff (\exists!\alpha)R(x,\alpha)\hspace*{4cm}
\end{multline*}
which is not very simple in this case because it requires constructive proofs
of Theorems~\ref{fact1} and~\ref{fact2}; but once thought out,
these proofs yield (practically) routinely the ``uniformities''
$\bu_1(\vep,x), \bu_2(\zeta)$ etc., and the precise statements
and rigorous proofs of these results and the Main Theorem. A
rigorous definition of ``constructive proof'' which justifies this
mode of reasoning is given in~\cite{ynmshanin2010}, but it is not
practically useful except for mathematicians who are familiar with
intuitionistic logic.

\bibliographystyle{ynmasl2}
\bibliography{sg}

\end{document}